\documentclass[12pt,russian,reqno]{amsart}

\usepackage[cp1251]{inputenc}
\usepackage[russian]{babel}
\usepackage{graphicx}

\usepackage{amsmath}
\usepackage{amsfonts}
\usepackage{amssymb}

\begin{document}

\noindent{УДК 514.743.4}

\begin{center}
\textbf{В.\,А.~Александров}$\,{}^{1,2}$
\end{center}

\medskip

{\small
\hfill${}^1\,$\textsf{Институт математики им. С.\,Л.~Соболева СО РАН}

\hfill\textsf{пр. Акад. Коптюга, 4, Новосибирск, 630090, Россия}

\hfill\textsf{E-mail: alex@math.nsc.ru}

\medskip	

\hfill${}^2\,$\textsf{Новосибирский государственный университет}

\hfill\textsf{ул. Пирогова, 2, Новосибирск, 630090, Россия}
}%endsmall

\newcounter{exer}
\addtocounter{exer}{1}
\addtocounter{footnote}{-1}

\bigskip

\bigskip

\begin{center}
\textbf{ПЕРВОЕ ЗНАКОМСТВО С ТЕНЗОРАМИ$^*$}\let\thefootnote\relax\footnote{${}^*$Работа выполнена при поддержке 
ФЦП <<Научные и научно-педагогические кадры инновационной России>> на 2009--2013 гг. 
(гос. контракт 02.740.11.0457) и Совета по грантам Президента РФ для поддержки ведущих
научных школ (НШ-6613.2010.1).\par{}\par
\texttt{ISSN 1818-7994. Вестник НГУ. Серия: Физика. 2012. Том 7, выпуск 1. Стр. 100--117.}}
\end{center}

\bigskip

{\small
\begin{center}
\begin{minipage}[t]{140mm}
Излагаются основы теории тензоров, в частности,
основы тензорной алгебры, значительное внимание уделелено
навыкам работы с символами Кронекера и Леви--Чивиты.
Приведено 60 упражнений для самостоятельного решения.
Статья адресована студентам младших курсов физических, 
математических и геолого-геофизических факультетов университетов.

\textit{Ключевые слова:} тензор, символ Кронекера, символ Леви--Чивиты.
\end{minipage}
\end{center}
}%endsmall

%Надо сегодня сказать лишь то, что уместно сегодня. 
%Прочее все отложить и сказать в подходящее время.
%--- Гораций (Квинт Гораций Флакк)
%Ut jam nunc dicat jam nunc debentia dici,
%Pleraque differat et praesens in tempus omittat
%[Надо сегодня сказать лишь то, что уместно сегодня,
%Прочее все отложить и сказать в подходящее время (лат.)].
%Ars Poetica, L. 43
%Let him say at first things that ought to be first said,
%and reserve for another time the gratest part of those
%that should  also have been said at first.
%The poetical works of ... Wentworth Dillon, earl of Roscommon, Page 137

\bigskip

\section{Введение}

Тензорное исчисление зарекомендовало себя как удобный 
математический инструмент, широко применяемыйв физике и геофизике. 
В том или ином виде тензоры довольно рано появляются в
современных университетских курсах. 
Например, тензоры инерции появляются в курсах механики и 
аналитической механики, а понятие 4-вектора --- в курсе 
электродинамики.
Но вот математический фундамент под эти понятия подводится 
значительно позже.

Наша задача --- дать студентам возможность познакомиться с 
тензорными методами на раннем этапе математического образования. 
Предпочтительно --- в конце второго семестра, когда уже пройдены 
основные разделы курса <<Линейная алгебра и геометрия>>. 
Из этого курса нам понадобятся самые общие сведения о системах 
координат, матрицах, определителях, скалярном и векторном 
произведениях. 
При необходимости читатель может освежить свои знания, 
обратившись к любому учебнику по линейной алгебре и геометрии
(см., например, к \cite{Ul07}).

Статья предназначена для первого знакомства с предметом
и не имеет целью дать полное изложение теории тензоров.
Большое внимание в ней уделелено изложению основ тензорной 
алгебры и навыкам работы с символами Кронекера и Леви--Чивиты.
Приведено 60 упражнений для самостоятельного решения.
Звёздочкой отмечены обязательные упражнения:
читателю, не научившемуся их решать, не следует думать,
что он освоил хотя бы основы теории тензоров.
Тем читателям, кто хочет более основательно разобраться в 
предмете, мы рекомендуем решить все упражнения.

Наконец, отметим, что ниже нам встретятся тензоры и другие 
объекты, названные в честь конкретных людей.
Это Леоп\'ольд Кр\'онекер (Leopold Kronecker, 1823--1891)~--- 
немецкий математик, специалист по алгебре и теории чисел и

Т\'уллио Л\'еви--Чив\'ита (Tullio Levi--Civita, 1873--1941)~--- 
итальянский математик и физик, один из создателей тензорного 
исчисления, много сделавший для его применения в теории 
относительности\footnote
{Напомним правило русского языка: 
мужские фамилии иностранного происхождения склоняются, 
а женские не склоняются.}.

\section{Индексные обозначения}

Система индексных обозначений составляет столь значительную 
часть тензорного исчисления, а технология манипулирования 
индексами столь востребована в физической литературе,
что знакомство с тензорами имеет смысл начать именно с неё.

\subsection*{Объекты.}
Совокупность трёх независимых переменных можно обозначить 
тремя различными буквами: например, $x, y, z$. 
Но во многих ситуациях более удобно обозначать переменные 
данной совокупности одной и той же буквой, различая их 
посредством индексов. 
В таком случае пишут, например, $x_1, x_2, x_3$, или в более 
компактной форме $x_r$, $r=1,2,3$. 
Записывая индекс $r$ внизу мы пока не вкладываем с это никакого 
особого смысла. 
В равной мере мы могли бы испольовать и верхний индекс, 
так что наша совокупность трёх переменных была бы записана 
в виде $x^1, x^2, x^3$, или $x^r$, $r=1,2,3$.
Нужно только понимать, что $x^r$ не означает возведение $x$ 
в $r$-ю степень, а уж если нам понадобится возвести переменную 
$x^r$ в степень $k$, то мы будем писать $(x^r)^k$.

Объекты, которые, подобно $x_r$ и $x^r$, зависят только от 
одного индекса, называют \textit{объектами первого порядка}, 
а отдельные буквы с индексами $x_1, x_2, x_3$ и 
$x^1, x^2, x^3$ называют \textit{компонентами} объекта. 
Объекты первого порядка, имеющие три составляющие, называют 
трёхмерными. 
Очевидно, имеется ровно два типа объектов первого порядка: 
те, у которых индекс внизу (например, $x_r$) и те, у которых 
индекс вверху (например, $x^r$).

Объекты, которые зависят от двух индексов, называются 
\textit{объектами второго порядка}. 
Поскольку индексы бывают верхние и нижние, то объекты второго 
порядка могут быть трёх типов: $x_{rs}$, $x_r^s$ и $x^{rs}$. 
Очевидно, каждый (трёхмерный) объект второго порядка имеет 9 
компонент.

Аналогично можно конструировать объекты любого порядка, 
содержащие любое количество индексов. 
Объект, зависящий от $m$ нижних индексов и от $n$ верхних 
называют \textit{объектом ранга} $(m,n)$.
Для того, чтобы придать этой системе обозначений завершённость,
объект $x$, не имеющий индексов, называют 
\textit{объектом нулевого порядка},
соответственно, он является объектом ранга $(0,0)$.

Для того, чтобы сделать формулы более компактными, 
в тензорном исчислении применяют следующие два соглашения 
относительно индексов:

1) \textit{повторяющийся малый латинский индекс означает 
сумирование от} 1 \textit{до} 3;

2) \textit{свободные} (\textit{неповторяющиеся}) 
\textit{малые латинские индексы пробегают значения от} 1 
\textit{до} 3.

Соглашение 1 часто называют <<\textit{правилом суммирования 
Эйнштейна}>>. 
В соответствии с ним развёрнутая запись выражения $a_rx^r$ 
имеет вид $a_1x^1+a_2x^2+a_3x^3$, а развёрнутая запись
выражения $a_{rs}x^rx^s$ имеет вид 
$$a_{11}(x^1)^2+a_{12}x^1x^2+a_{13}x^1x^3+a_{21}x^2x^1+a_{22}(x^2)^2+a_{23}x^2
x^3+a_{31}x^3x^1+a_{32}x^3x^2+a_{33}(x^3)^2.$$
Поскольку по повторяющемуся индексу подразумевается сумирование 
от 1 до 3, то, заменив его любой другой буквой, мы не изменим 
значения рассматриваемого выражения (тем самым, например, 
$a_rx^r \equiv a_sx^s$ и $a_{rs}x^rx^s\equiv a_{mn}x^mx^n$).
По этой причине повторяющийся индекс часто называют 
\textit{немым}. 
Сокращённая запись, в которой один и тот же немой индекс 
появляется более двух раз (например $x_{rrr}$), считается 
запрещённой.
Соглашение 2 освобождает нас от необходимости после каждой 
формулы указывать пределы изменениям свободных индексов.

\medskip

\noindent\textit{Упражнения}

\noindent\textbf{\theexer.} Выпишите в развёрнутом виде 
систему линейных равенств, задаваемую выражением
$a_{rs}x^s=b_r$.
\addtocounter{exer}{1}

\noindent\textbf{\theexer.} Сколько членов содержится в сумме 
$a_{rst}x^ry^sz^t$?
\addtocounter{exer}{1}

\noindent\textbf{\theexer.} Покажите, что 
$x_{\hphantom{r}s}^{rs}$ есть объект первого порядка и 
выпишите все его компоненты.
\addtocounter{exer}{1}

\medskip

\subsection*{Операции над объектами.}
На множестве объектов со многими индексами имеются три основные 
операции, называемые сложением, умножением и свёрткой. 
Дадим соответствующие определения.

Если нам даны два объекта одного и того же ранга, то, сложив 
каждую компоненту первого объекта с соответствующей компонентой 
второго, мы, очевидно, придём к объекту того же ранга, что и 
слагаемые. 
Он называется \textit{суммой} этих обектов.
Например: если $x_{st}^r$ и $y_{st}^r$ --- два объекта третьего 
порядка, то объект $z_{st}^r$, определённый равенством 
$z_{st}^r=x_{st}^r+y_{st}^r$, есть сумма объектов $x_{st}^r$ и 
$y_{st}^r$.

Если взять два объекта произвольных рангов и умножить каждую 
компоненту первого объекта на каждую компоненту второго, то 
результирующий объект называют \textit{произведением}
исходных двух объектов. 
Очевидно, порядок произведения равен сумме порядков исходных 
объектов. 
Например: если $x_{st}^r$ --- объект третьего порядка и 
$y^{mn}$ --- объект второго порядка, то объект $z_{st}^{rmn}$, 
определённый равенством $z_{st}^{rmn}=x_{st}^ry^{mn}$, 
есть произведение объектов $x_{st}^r$ и $y^{mn}$.

Построение свёртки поясним на примере. Если нам дан объект 
$x_{kst}^{rp}$ ранга $(3,2)$, то мы можем положить $k$ равным 
$p$ и получим объект $x_{pst}^{rp}$. Поскольку теперь
$p$ является повторяющимся индексом, то, в соответствии с 
соглашением 1, по нему нужно просуммировать от 1 до 3. 
Следовательно, полученный в результате объект
$$x_{pst}^{rp}\equiv x_{1st}^{r1}+x_{2st}^{r2}+x_{3st}^{r3}$$
имеет ранг $(2,1)$. 
Он называется \textit{свёрткой} объекта $x_{kst}^{rp}$ 
по индексам $k$ и $p$.
Очевидно, операция свёртки может быть повторена несколько раз.
При этом каждый раз в паре сворачиваемых индексов один должен 
быть верхним, а другой --- нижним.

Объект называется \textit{симметричным} относительно двух 
своих нижних индексов, если его компоненты не изменяются 
при перемене мест этих двух индексов. 
Например, объект $x_{rs}$ является симметричным, если 
$x_{rs}=x_{sr}$.
Объект называется \textit{абсолютно симметричным} относительно 
нижних индексов, если при перемене мест любых двух из них 
компоненты не меняются. 
Например, абсолютно симметричный объект $x_{rst}$ удовлетворяет 
равенствам $$x_{rst}=x_{rts}=x_{str}=x_{srt}=x_{trs}=x_{tsr}.$$

Объект называется \textit{антисимметричным} относительно двух 
своих нижних индексов, если перемена мест этих двух индексов 
изменяет знак компоненты. 
Например, объект $x_{rs}$ является антисимметричным, если 
$x_{rs}=-x_{sr}$.
Объект называется \textit{абсолютно антисимметричным} 
относительно нижних индексов, если он антисимметричен 
относительно всех пар нижних индексов. 
Например, абсолютно симметричный объект $x_{rst}$ 
удовлетворяет равенствам 
$$x_{rst}=-x_{rts}=x_{str}=-x_{srt}=x_{trs}=-x_{tsr}.$$

Всё сказанное выше о симметрии и антисимметрии объектов 
в равной степени применимо и к верхним индексам.

\medskip

\noindent\textit{Упражнения}

\noindent\textbf{\theexer.} Сколько различных компонент имеет 
абсолютно симметричный объект $x_{rst}$?
\addtocounter{exer}{1}

\noindent\textbf{\theexer$^\ast$.} Покажите, что абсолютно 
антисимметричный (трёхмерный) объект $y_{rst}$ имеет самое 
б\'ольшее шесть отличных от нуля компонет, причём все они 
одинаковы по абсолютной величине.
\addtocounter{exer}{1}

\noindent\textbf{\theexer.} Докажите, что если объект $a_{rs}$ 
антисимметричен, то $a_{rs}x^rx^s=0$; и обратно, если это 
уравнение верно для всех объектов $x^r$, то объект $a_{rs}$ 
антисимметричен.
\addtocounter{exer}{1}

\subsection*{Символы Кронекера и Леви--Чивиты.}
Определим три объекта $\delta_{rs}$,  $\delta^{rs}$ и 
$\delta_r^s$ условием: компонента с индексами $r$ и $s$ равна 1, 
если $r=s$ и равна 0, если $r\neq s$. 
Каждый их этих объектов называется 
\textit{символом Кр{\'о}некера}. 
Очевидно, каждый из объектов $\delta_{rs}$ и $\delta^{rs}$ 
является симметричным.

Определим объекты $e_{rst}$ и  $e^{rst}$ условием: 
компонента с индексами $r$, $s$ и $t$ равна 0, 
если среди индексов есть хотя бы два одинаковых; 
равна 1, если тройка чисел $r$, $s$ и $t$ является чётной 
перестановкой чисел 1, 2 и 3; 
и равна $-1$, если тройка чисел $r$, $s$ и $t$ является 
нечётной перестановкой чисел 1, 2 и 3. 
Каждый их этих объектов называется 
\textit{символом Л{\'е}ви--Чив{\'и}ты}. 
Очевидно, каждый из объектов $e_{rst}$ и  $e^{rst}$ является 
абсолютно антисимметричным. 

\medskip

\noindent\textit{Упражнения}

\noindent\textbf{\theexer.} Докажите равенство $\delta_r^r=3$.
\addtocounter{exer}{1}

\noindent\textbf{\theexer$^\ast$.} Докажите, что для 
произвольного объекта $x^r$ имеет место равенство 
$\delta_s^r x^s=x^r$.
\addtocounter{exer}{1}

\noindent\textbf{\theexer$^\ast$.} Докажите, что если 
$x_{rst}$ --- произвольный абсолютно антисимметричный объект, 
то $x_{rst}=x_{123}e_{rst}$.
\addtocounter{exer}{1}

\noindent\textbf{\theexer.} Докажите формулу 
$e_{rst}=e^{rst}=(s-r)(t-r)(t-s)/2$.
\addtocounter{exer}{1}

\subsection*{Индексные обозначения и теория определителей.}
Покажем как введённые выше индексные обозначения позволяют 
по-новому взглянуть на известные вам
факты теории определителей.

Определитель, элементами которого являются компоненты 
трёхмерного объекта $x_s^r$, будем обозначать через
$$
\det x_s^r = 
\det\begin{pmatrix}
x_1^1 & x_2^1 & x_3^1 \\
x_1^2 & x_2^2 & x_3^2 \\
x_1^3 & x_2^3 & x_3^3 
\end{pmatrix}.
$$
(Тем самым мы подразумеваем, что верхний индекс объекта 
$x_s^r$ обозначает строку, а нижний --- столбец соответствующей 
матрицы.)

Из линейной алгебры известно, что определителем такой матрицы 
называется сумма произведений
$x_1^rx_2^sx_3^t$ элементов матрицы, взятых по одному 
из каждой строки и каждого столбца,
причём каждое произведение берётся со знаком <<$+$>>, 
если тройка чисел $r$, $s$ и $t$ является чётной перестановкой 
чисел 1, 2 и 3 и со знаком <<$-$>>, если тройка чисел 
$r$, $s$ и $t$ является нечётной перестановкой чисел 1, 2 и 3.
Следовательно, принимая во внимание определение символа 
Леви--Чивиты $e_{rst}$ и правило суммирования Эйнштейна, 
можем написать 
$$\det x_s^r = e_{rst} x_1^rx_2^sx_3^t.\eqno(1)$$

Отметим, что из этого равенства сразу вытекает хорошо известное 
свойство определителя:
\textit{При перемене мест любых двух столбцов определитель 
меняет знак.} 

Пользуясь тем, что значение суммы не изменится, если мы 
используем другую букву для обозначения 
немой переменной, и пользуясь тем, что символ Леви--Чивиты 
$e_{rst}$ является абсолютно антисимметричным, можем написать
$$e_{rst} x_m^rx_n^sx_p^t=e_{tsr} x_m^tx_n^sx_p^r
=e_{tsr} x_n^sx_p^rx_m^t=-e_{rst} x_n^sx_p^rx_m^t.$$
Значит, объект $e_{rst} x_m^rx_n^sx_p^t$
меняет знак при перестановке индексов $r$ и $t$. 
Ясно, что мы получим тот же результат, если
поменяем местами любые два из индексов $r$, $s$, $t$. 
Тем самым мы убеждаемся, что объект $e_{rst} x_m^rx_n^sx_p^t$ 
является абсолютно антисимметричным.

Наконец, используя упражнение 9 и формулу (1)  получим
$$
e_{rst} x_m^rx_n^sx_p^t
=e_{mnp} e_{rst} x_1^rx_2^sx_3^t
=e_{mnp} \det x_j^k.\eqno(2)
$$

Теперь всё готово для доказательства хорошо известной теоремы 
об умножении определителей, которую мы сформулируем так: 
\textit{Пусть даны три объекта $x_s^r$, $y_s^r$ и $z_s^r$, 
причём $z_s^r=x_m^ry_s^m$} (т.\,е. если мыслить эти объекты 
как матрицы
$$
\begin{pmatrix}
x_1^1 & x_2^1 & x_3^1 \\
x_1^2 & x_2^2 & x_3^2 \\
x_1^3 & x_2^3 & x_3^3 
\end{pmatrix},\quad
\begin{pmatrix}
y_1^1 & y_2^1 & y_3^1 \\
y_1^2 & y_2^2 & y_3^2 \\
y_1^3 & y_2^3 & y_3^3 
\end{pmatrix}\quad\mbox{и}\quad
\begin{pmatrix}
z_1^1 & z_2^1 & z_3^1 \\
z_1^2 & z_2^2 & z_3^2 \\
z_1^3 & z_2^3 & z_3^3 
\end{pmatrix},
$$
то матрица объекта $z_s^r$ является произведением матриц 
объектов $x_s^r$ и $y_s^r$). \textit{Тогда}
$\det z_s^r = (\det x_s^r)(\det y_s^r).$

Доказательство получается прямым вычислением
\begin{align*}
(\det x_s^r)(\det y_s^r)&=
(\det x_s^r)(e_{mnp} y_1^my_2^ny_3^p)=
\bigl((\det x_s^r)e_{mnp}\bigr)(y_1^my_2^ny_3^p)=\\
&=(e_{rst} x_m^rx_n^sx_p^t)(y_1^my_2^ny_3^p)=
e_{rst} (x_m^ry_1^m)(x_n^sy_2^n)(x_p^ty_3^p)
=e_{rst} z_1^rz_2^sz_3^t=\det z_s^r,
\end{align*}
в котором первое и последнее равенства написаны в силу 
формулы (1), третье --- в силу формулы (2), а второе и 
четвёртое отражают перегруппировку слагаемых.

\medskip

\noindent\textit{Упражнения}

\noindent\textbf{\theexer$^\ast$.} Докажите равенство 
$\det \delta_s^r=1$.
\addtocounter{exer}{1}

\noindent\textbf{\theexer$^\ast$.} Докажите, что если объект 
$x_s^r$ обладает свойством 
$x_m^rx_s^m=\delta_s^r$, то $\det x_s^r=\pm 1$. 
Сформулируйте аналогичное свойство матриц,
известное вам из курса линейной алгебры.
\addtocounter{exer}{1}

\noindent\textbf{\theexer.} Докажите равенство 
$\det x_s^r = e^{rst} x^1_rx^2_sx^3_t$,
аналогичное равенству (1).
Выведите из него хорошо известное свойство определителя:
\textit{При перемене мест любых двух строк определитель 
меняет знак.}
\addtocounter{exer}{1}

\noindent\textbf{\theexer.} Докажите равенство 
$e^{rst} x^m_rx^n_sx^p_t
=e^{mnp} \det x_j^k,$ аналогичное равенству (2).
\addtocounter{exer}{1}

\noindent\textbf{\theexer$^\ast$.} Докажите равенство 
$$
e_{mnp}e^{rst} = 
\det\begin{pmatrix}
\delta_m^r & \delta_n^r& \delta_p^r\\
\delta_m^s & \delta_n^s& \delta_p^s\\
\delta_m^t & \delta_n^t& \delta_p^t\\
\end{pmatrix}.\eqno(3)
$$
[\texttt{Указание:} Убедитесь, что левая и правая части равенства (3) 
являются абсолютно антисимметричными объектами по всем нижним 
(и по всем верхним) индексам; выведите отсюда, что (3) 
достаточно проверить при каком-нибудь
одном наборе индексов, при котором левая часть не равна нулю.]
\addtocounter{exer}{1}

\noindent\textbf{\theexer$^\ast$.} Докажите равенство 
$e_{mnp}e^{rsp}=\delta_m^r\delta_n^s-\delta_m^s\delta_n^r.$\newline
{\hphantom{a}}\hfill [\texttt{Указание:} Воспользуйтесь равенством (3).]
\addtocounter{exer}{1}

\noindent\textbf{\theexer.} Докажите равенство 
$e_{mnp}e^{rnp}=2\delta_m^r.$\newline
{\hphantom{a}}\hfill [\texttt{Указание:} Воспользуйтесь предыдущим упражнением.]
\addtocounter{exer}{1}

\noindent\textbf{\theexer.} Докажите равенство 
$e_{mnp}e^{mnp}=6.$\newline
{\hphantom{a}}\hfill [\texttt{Указание:} Воспользуйтесь предыдущим упражнением.]
\addtocounter{exer}{1}

\section{Основы тензорной алгебры}

Несколько опережая события можно сказать, что тензор есть 
частный случай объекта, рассмотренного в предыдущем разделе, 
компоненты которого при замене координат
преобразуются по некоторому определённому закону.
Самое важное при этом --- понять зачем вообще нужно следить 
за тем, как преобразуются компоненты при замене координат 
и почему они должны преобразовываться именно
по тому закону, который мы укажем в своё время в определении.

Чтобы разобраться в этом, обратимся к примерам.

\subsection*{Контравариантный вектор.}
Начнём с того, что <<вектор --- это направленный отрезок>>. 
Наверное, ещё в школе вы узнали, что сила --- это вектор
и уже тогда не раз представляли её себе в виде <<направленного 
отрезка>>, приложенного к какой-то точке тела, и складывали 
такие направленные отрезки по правилу параллелограмма.
Другой полезный для нас сейчас пример --- это радиус-вектор 
точки, т.\,е. <<направленный отрезок>>, идущий из начала 
координат в данную точку.

У такого геометрического подхода к векторам есть существенный 
недостаток: <<направленный отрезок>> нельзя подставить в 
уравнение. 
Вы уже знаете, что в уравнения входят координаты векторов. 
Но координаты вектора вычисляются относительно какой-то
системы координат. 
Однако сами по себе системы координат в природе не встречаются. 
По мере надобности их вводят в рассмотрение исследователи.
В частности, один исследователь может использовать одну 
систему координат, а другой --- другую. 
И вот тут-то и начинается самое интересное: вектор (скажем, 
вектор силы или радиус-вектор точки) даны нам <<реально и 
объективно>>, а координаты этого вектора у каждого исследователя 
свои (скажем потому, что каждый использует свою, любимую
систему координат). 
Как же этим исследователям договориться и понять, что они 
работают с одним и тем же <<реальным>> вектором? 
Ответ на этот вопрос даёт <<закон преобразования>> координат 
вектора, к выводу которого мы и переходим.

Пусть в трёхмерном евклидовом пространстве заданы две системы 
координат, имеющие общее начало. 
Одну из них будем называть <<старой>> системой координат 
и в ней будем обозначать координаты произвольной точки 
пространства через $x^1$, $x^2$, $x^3$. 
Другую систему координат будем называть <<новой>> и
в ней будем обозначать координаты произвольной точки через
$\overline{x}^1$, $\overline{x}^2$, $\overline{x}^3$.
Особо отметим, что здесь и далее черта не означает 
комплексного сопряжения.

Будем считать, что <<старые>> координаты $x^1$, $x^2$, $x^3$
преобразуются в <<новые>> $\overline{x}^1$, $\overline{x}^2$, 
$\overline{x}^3$
с помощью линейного невырожденного преобразования
\begin{align*}
\overline{x}^1&= c_1^1x^1+c_2^1x^2+c_3^1x^3,\\
\overline{x}^2&= c_1^2x^1+c_2^2x^2+c_3^2x^3,\\
\overline{x}^3&= c_1^3x^1+c_2^3x^2+c_3^3x^3,
\end{align*}
где $c_r^s$ --- некоторые постоянные.
Используя правило суммирования Эйнштейна мы можем записать 
эту систему уравнений в виде
$$
\overline{x}^r= c_s^rx^s.\eqno(4)
$$

Таким образом мы видим, что существуют <<реальные>> объекты 
(например, --- радиус-вектор точки), компоненты которых 
преобразуются по закону (4).
Любой объект $x^r$, компоненты которого преобразуются по 
закону (4), называется \textit{контравариантным тензором 
первого ранга}, или, иначе, \textit{контравариантным вектором}.

\subsection*{Ковариантный вектор.}
Зададимся вопросом о том все ли объекты первого порядка
таковы, что их компоненты преобразуются по закону (4)?

Прежде всего заметим, что поскольку преобразование (4) 
предположено невырожденным, то
$\det c_r^s\neq 0$ и можно разрешить систему уравнений (4) 
относительно переменных $x^1$, $x^2$, $x^3$:
$$
x^r= \gamma_s^r\overline{x}^s,\eqno(5)
$$
где $\gamma_s^r$ --- некоторые постоянные. 
(На самом деле они вполне определённые:
матрица $\gamma_s^r$ является в точности обратной 
к матрице $c_s^r$; см. упражнение 20.)

Предположим теперь, что в пространстве задан линейный 
функционал, который задаётся в <<старой>> 
системе координат формулой $a_rx^r$. 
Мы предполагаем этот функционал инвариатным,
т.\,е. таким, что его значение в конкретной точке не 
зависит от выбора какой-либо системы координат. 
Наша задача состоит в том, чтобы найти компоненты 
$\overline{a}_r$ той линейной функции
$\overline{a}_r\overline{x}^r$, которая задёт этот же 
самый функционал в <<новой>> системе координат.

Поскольку функционал инвариантен, то 
$a_rx^r= \overline{a}_r\overline{x}^r$. 
Воспользовавшись (5), отсюда получим 
$a_r\gamma_s^r\overline{x}^s= \overline{a}_r\overline{x}^r$,
или, переобозначив немые индексы, 
$$(a_s\gamma_r^s- \overline{a}_r)\overline{x}^r=0.\eqno(6)$$
Формула (6) должна иметь место при любом выборе численных 
значений координат
$\overline{x}^1$, $\overline{x}^2$, $\overline{x}^3$.
Воспользуемся этим обстоятельством: фиксировав $r$, 
положим $\overline{x}^r=1$, а все 
остальные координаты $\overline{x}^p$, $p\neq r$, 
положим равными нулю.
Для такого набора координат непосредственно из (6) получим
$$
\overline{a}_r=\gamma_r^sa_s.\eqno(7)
$$

Следовательно, существуют <<реальные>> объекты первого порядка
(например, линейные функционалы), компоненты
которых преобразуются не по закону (4), а по закону (7). 
Любой объект $a_r$, компоненты которого преобразуются по 
закону (7), называется \textit{ковариантным тензором первого 
ранга}, или, иначе, \textit{ковариантным вектором}.

Таким образом у нас есть два типа тензоров первого порядка, и мы 
условимся различать их с помощью положения индекса:
\textit{если тензор контравариантен, то мы используем верхний 
индекс\/}; \textit{если же он ковариантен, то нижний}.

\medskip

\noindent\textit{Упражнения}

\noindent\textbf{\theexer$^\ast$.} Докажите формулу 
$x_r=c_r^s\overline{x}_s$.
\addtocounter{exer}{1}

\noindent\textbf{\theexer$^\ast$.} Пусть объекты  
$c_r^s$ и $\gamma_r^s$ взяты из формул (4) и (5). 
Докажите формулы $\gamma_r^sc_s^t=\delta_r^t$ и
$\gamma_r^sc_t^r=\delta_t^s$, означающие, что матрица 
$\gamma_s^r$ является обратной к матрице $c_s^r$.
\addtocounter{exer}{1}

\subsection*{Тензор ранга (2,0).}
Теперь перейдём к рассмотрению известных вам из курса 
линейной алгебры инвариантных объектов с двумя индексами. 
В частности, вы изучали кваратичные формы
и законы их преобразования к новым переменным. 
Вот с них мы и начнём, но 
теперь будем использовать индексные обозначения. 

Итак, пусть имеется инвариантная квадратичная форма, 
заданая в <<старых>> координатах $x^r$ формулой $a_{rs}x^rx^s$, 
а в <<новых>> координатах $\overline{x}^s$ ~--- формулой
$\overline{a}_{rs}\overline{x}^r\overline{x}^s$.
Изначально будем предполагать, что компоненты формы 
симметричны, т.\,е. 
$$a_{rs}=a_{sr} \quad\mbox{и}\quad \overline{a}_{rs}=\overline{a}_{rs}.\eqno(8)$$
Отметим, что предположение (8) не ведёт к ограничению общности, 
поскольку если, например, $a_{12}\neq a_{21}$, то, заменив 
коэффициенты $a_{12}$ и $a_{21}$ новыми по формулам 
$$\widetilde{a}_{12}=\widetilde{a}_{21}=(a_{12}+a_{21})/2,$$
мы не изменим нашу форму как однородную квадратичную функцию, 
сопоставлющую точке пространства вещественное число, 
но сделаем коэффициенты с индексами {\tiny 12} и 
{\tiny 21} одинаковыми (а потом повторим эту процедуру 
для всех остальных пар индексов). 

Итак, поскольку наша форму инвариантна, то 
$a_{rs}x^rx^s=\overline{a}_{rs}\overline{x}^r\overline{x}^s$. 
Используя формулу (5), можем переписать её в виде 
$$(a_{rs}\gamma_m^r\gamma_n^s-\overline{a}_{mn})x^mx^n=0.\eqno(9)$$
Поскольку равенство (9) должно выполняться при всех значениях координат 
$\overline{x}^1$, $\overline{x}^2$, $\overline{x}^3$,
то, фиксировав два разных индекса $m$ и $n$, вычислим значение 
нашей формы в той точке пространства, у которой $m$-я и $n$-я 
<<новые>> координаты равны единице, а третья координата,
имеющая индекс, отличный от $m$ и $n$, равна нулю. 
Получим 
$$\gamma_m^r\gamma_n^sa_{rs}+\gamma_m^s\gamma_n^ra_{sr}
-\overline{a}_{mn}-\overline{a}_{nm}=0,$$
или, используя предположения (8) и выполняя очевидные 
преобразования, 
$$
\overline{a}_{mn}=\gamma_m^r\gamma_n^sa_{rs}.\eqno(10)
$$
Это и есть искомый закон преобразования коэффициентов 
квадратичной формы.
Любой объект $a_{rs}$, компоненты которого преобразуются 
по закону (10), называется \textit{тензором ранга} (2,0) 
или \textit{ковариантным тензором второго порядка}.

\subsection*{Преобразование векторов базиса.}
Прежде чем двигаться дальше, выясним как преобразуются векторы 
базиса при замене координат в пространстве.
Выше мы очень кратко сказали, что <<старые>> координаты 
$x^1$, $x^2$, $x^3$ преобразуются в <<новые>> 
$\overline{x}^1$, $\overline{x}^2$, $\overline{x}^3$
с помощью линейного невырожденного преобразования (4).
Но из линейной алгебры вы знаете, что задание системы 
координат в трёхмерном пространстве осуществляется заданием 
тройки некомпланарных векторов, называемых 
\textit{репером} \cite{Ul07}\footnote{Эта
книга опубликована также в Интернете:
\texttt{http://www.phys.nsu.ru/ok03/manuals\_2007-2008.html}}.
Значит, при написании формулы (4) мы подразумевали, что у 
нас есть <<старая>> система координат  $x^1$, $x^2$, $x^3$, 
соответствующая реперу
$$\boldsymbol{e}_1, \boldsymbol{e}_2, \boldsymbol{e}_3,\eqno(11)$$
и есть <<новая>> система координат  
$\overline{x}^1$, $\overline{x}^2$, $\overline{x}^3$, 
соответствующая реперу
$$\overline{\boldsymbol{e}}_1, \overline{\boldsymbol{e}}_2, \overline{\boldsymbol{e}}_3.\eqno(12)$$

Мы знаем, что координаты  $x^1$, $x^2$, $x^3$ и 
$\overline{x}^1$, $\overline{x}^2$, $\overline{x}^3$ 
связаны между собой соотношениями (4) и (5) и хотим найти 
соотношения, связывающие векторы реперов (11) и (12).

Поскольку тройка векторов (12) образует базис в пространстве, 
то всякий вектор представляется
в виде их линейной комбинации; в частности,
\begin{align*}
\boldsymbol{e}_1&=\beta_1^1\overline{\boldsymbol{e}}_1+\beta_1^2\overline{\boldsymbol{e}}_2+
\beta_1^3\overline{\boldsymbol{e}}_3,\\
\boldsymbol{e}_2&=\beta_2^1\overline{\boldsymbol{e}}_1+\beta_2^2\overline{\boldsymbol{e}}_2+
\beta_2^3\overline{\boldsymbol{e}}_3,\\
\boldsymbol{e}_3&=\beta_3^1\overline{\boldsymbol{e}}_1+\beta_3^2\overline{\boldsymbol{e}}_2+
\beta_3^3\overline{\boldsymbol{e}}_3,
\end{align*}
где $\beta_r^s$ --- некоторые постоянные, которые мы намерены 
выразить через постоянные из формул (4) и (5). 
Пользуясь правилом суммирования Эйнштейна будем записывать 
последние формулы короче: 
$\boldsymbol{e}_r= \beta_r^s\overline{\boldsymbol{e}}_s$.

Разлагая радиус-вектор произвольной точки пространства один 
раз по базису (11), а второй раз по базису (12), получим
$x^r\boldsymbol{e}_r=\overline{x}^r\overline{\boldsymbol{e}}_r$.
Подставив сюда $\boldsymbol{e}_r$ из предыдущей формулы и 
$\overline{x}^r$ из формулы (4), будем иметь 
$x^r\beta_r^s\overline{\boldsymbol{e}}_s=
c_s^rx^s\overline{\boldsymbol{e}}_r$, или,
поменяв ролями индексы $r$ и $s$ в правой части,
$x^r\beta_r^s\overline{\boldsymbol{e}}_s=
c_r^sx^r\overline{\boldsymbol{e}}_s$.
Поскольку любой вектор может быть разложен по базису (12) 
единственным образом, то из последней формулы заключаем, что 
$x^r\beta_r^s=c_r^sx^r$.
А поскольку $x^r$ произволен\footnote{Так коротко ссылаются
на приём, которым мы уже не раз пользовались выше и который 
состоит в том, что, фиксировав индекс $r$, нужно применить 
последнюю формулу к такой точке пространства, у которой 
$r$-я координата равна единице, а все остальные --- нулю; 
каждая из сумм $x^r\beta_r^s=c_r^sx^r$ сведётся к одному 
слагаемому и мы получим $\beta_r^s=c_r^s$.}, 
то отсюда сразу выводим, что $\beta_r^s=c_r^s$.
Это и есть искомое соотношение, которое можно записать ещё и так:
$$\boldsymbol{e}_r= c_r^s\overline{\boldsymbol{e}}_s.\eqno(13)$$

\medskip

\noindent\textit{Упражнение}

\noindent\textbf{\theexer$^\ast$.} Докажите, что если <<старые>> 
и <<новые>> координаты связаны между собой соотношениями (4) 
(или, что то же самое, соотношениями (5)), то
$\overline{\boldsymbol{e}}_r= \gamma_r^s\boldsymbol{e}_s.$
\addtocounter{exer}{1}

\subsection*{Тензор ранга (1,1).}
Теперь мы можем перейти к рассмотрению ещё одного известного 
вам из курса линейной алгебры инвариантного объекта с двумя 
индексами. 
Речь пойдёт о линейных отображениях. 
Как вы знаете, линейное отображение (или оператор) представляет 
собой правило, сопоставляющее каждому вектору трёхмерного 
пространства другой вектор этого же пространства (конечно, 
это правило таково, что соблюдается условие линейности). 
Инвариантность оператора состоит в том, что это правило 
не зависит от выбора системы координат. 

Вы знаете, что если в пространстве введена система координат, 
то оператору однозначно соответствует \textit{матрица оператора} 
$a_s^r$, так что компоненты $x^1$, $x^2$, $x^3$ любого вектора 
связаны с компонентами $y^1$, $y^2$, $y^3$ его образа под 
действием нашего оператора посредством формулы $y^r=a_s^rx^s$. 
Точно так же в <<новой>> системе координат тому же самому 
оператору соответствует матрица $\overline{a}_s^r$, 
так что компоненты $\overline{x}^1$, $\overline{x}^2$, 
$\overline{x}^3$ любого вектора связаны с компонентами 
$\overline{y}^1$, $\overline{y}^2$, $\overline{y}^3$ 
его образа под действием нашего оператора посредством формулы 
$\overline{y}^r=\overline{a}_s^r\overline{x}^s$. 
Наша задача состоит в том, чтобы выяснить, каким образом 
связаны между собой компоненты объектов $a_s^r$ и 
$\overline{a}_s^r$, если <<новые>> координаты
$\overline{x}^r$ выражаются через <<старые>> $x^r$ формулами (4).

Итак, поскольку линейное отображение инвариантно, то 
$y^r\boldsymbol{e}_r=\overline{y}^r\overline{\boldsymbol{e}}_r$, 
или, что то же самое,
$a_s^rx^s\boldsymbol{e}_r
=\overline{a}_s^r\overline{x}^s\overline{\boldsymbol{e}}_r$.
Подставив сюда выражения для $x^s$ и $\boldsymbol{e}_r$ 
из формул (5) и (14) соответственно, получим
$$a_s^r(\gamma_t^s\overline{x}^t)
(c_r^m\overline{\boldsymbol{e}}_m)=
\overline{a}_s^r\overline{x}^s\overline{\boldsymbol{e}}_r,$$
или, переобозначив немые индексы в левой части последнего равенства,
$$\gamma_s^tc_m^ra_t^m\overline{x}^s\overline{\boldsymbol{e}}_r=
\overline{a}_s^r\overline{x}^s\overline{\boldsymbol{e}}_r.$$
Приравнивая коэффициенты при одинаковых векторах базиса 
$\overline{\boldsymbol{e}}_r$ и пользуясь, как и выше, 
произвольностью координат $\overline{x}^s$, получим
$$\overline{a}_s^r=\gamma_s^tc_m^ra_t^m.\eqno(14)$$
Это и есть искомый закон преобразования матрицы линейного 
отображения.
Любой объект $a_s^r$, компоненты которого преобразуются 
по закону (14), называется \textit{тензором ранга} (1,1) 
или \textit{смешанным тензором второго порядка}.

\subsection*{Тензоры произвольного ранга.}
В самых общих чертах можно сказать, что тензор есть 
многомерный аналог вектора.
Для читателя, разобравшего приведённые выше примеры  
и проследившего вывод законов преобразования (4), (7), (10) и 
(14), следующее определение покажется совершенно естественным:

\textit{Говорят, что задан тензор ранга} $(m,n)$,
\textit{если в каждой системе координат трёхмерного пространства
задан объект $a_{i_1\dots i_m}^{j_1\dots j_n}$ ранга $(m,n)$,
т.\,е. объект, зависящий от $m$ нижних индексов и от $n$ 
верхних, который при задаваемом формулами} (4) 
\textit{переходе от <<старой>> сиcтемы координат к <<новой>>
преобразуется по правилу}
$$
\overline{a}_{i_1\dots i_m}^{j_1\dots j_n}=
\gamma_{i_1}^{r_1}\dots \gamma_{i_m}^{r_m}c_{s_1}^{j_1}\dots c_{s_n}^{j_n}
a_{r_1\dots r_m}^{s_1\dots s_n}.\eqno(15)$$

При этом компоненты объекта $a_{i_1\dots i_m}^{j_1\dots j_n}$
называются \textit{компонентами тензора} (если нужно, при этом 
указывают о какой системе координат идёт речь).
Нижние индексы компонент тензора называются 
\textit{ковариантными}, а верхние --- \textit{контравариантными}. 

Мнемоническое правило для запоминания тензорного закона 
преобразования (15) состоит в том, что ковариантные индексы 
преобразуются с помощью той же матрицы~$\gamma_r^s$, 
посредством которой преобразуются ковариантные векторы 
(см. формулу~(7)), а контавариантные индексы --- с помощью 
матрицы~$c_r^s$, посредством которой преобразуются
контравариантные векторы (см. формулу~(4)). 
Таким образом, на практике достаточно помнить формулы~(4) и~(7).

Из определения ясно, что тензоры могут быть любых рангов. 
Особо подчеркнём, что бывают и тензоры ранга $(0,0)$.
Они называются \textit{скалярами} или \textit{скалярными 
величинами}.
В соответствии с определением, чтобы задать скалярную 
величину, нужно в каждой системе координат задать объект $a$, 
не имеющий ни одного индекса, который, к тому же, не изменяется 
при заменах координат (ведь в этом случае $m=n=0$ и формула (15) 
принимает вид $\overline{a}=a$).

Рассмотрим подробнее один (но очень важный) пример скаляра.
Речь пойдёт о \textit{следе} $\mbox{tr\,} a_r^s$ тензора 
ранга (1,1) $a_r^s$.
В каждой системе координат мы определим его как свёртку 
соответствующего  объекта $a_r^s$ по индексам $r$ и $s$, 
т.\,е. по определению положим $\mbox{tr\,} a_r^s=a_r^r$. 
Убедимся, что так определённая величина действительно не 
зависит от выбора системы координат (и тем самым является 
скаляром). 
Это следует из прямого вычисления, основанного на формуле~(14) 
и упражнении~20:
$$\mbox{tr\,} \overline{a}_r^s=\overline{a}_r^r=\gamma_r^tc_p^ra_t^p=\delta_p^t a_t^p=
a_t^t=\mbox{tr\,} a_r^s.$$

Поскольку в ранее мы интерпретировали объект ранга (1,1)
как матрицу линейного отображения, то последннее вычисление 
является не чем иным, как  ещё одним доказательством уже 
известного вам из линейной алгебры утверждения, что
\textit{след линейного отображения не зависит от выбора системы 
координат}, т.\,е. инвариантен.

\medskip

\noindent\textit{Упражнения}

\noindent\textbf{\theexer.} Напишите законы преобразования 
тензоров $x^{rs}$ и $x_{st}^r$.
\addtocounter{exer}{1}

\noindent\textbf{\theexer.} Покажите, что если имеет место 
соотношение $x_{st}^r=y_s^rz_t$, связывающее компоненты 
тензоров $x_{st}^r$, $y_s^r$, $z_t$ в некоторой системе 
координат, то то же самое соотношение между компонентами 
имеет место в любой другой системе координат.
\addtocounter{exer}{1}

\noindent\textbf{\theexer$^\ast$.} Докажите, что если компоненты 
тензора с некоторой системе координат образуют симметричный 
(или антисимметричный) объект, то в любой другой системе 
координат компоненты этого тензора также образуют симметричный 
(или антисимметричный) объект. Такой тензор называют 
\textit{симметричным} (или \textit{антисимметричным}).
\addtocounter{exer}{1}

\subsection*{Простейшие свойства тензоров.} Пожалуй, самым 
важным является следующее почти очевидное свойство: 
\textit{если все компоненты тензора равны нулю в какой-то 
системе координат, то все его компоненты равны нулю в любой 
системе координат.}

Вместе с тем, указанное свойство является очень простым, 
ведь его доказательство немедленно вытекает из формулы~(15).

Другим важным (и опять-таки, очень простым) свойством тензоров 
является их \textit{произвольность}.
Под этим подразумевают, что \textit{какую бы систему координат 
и какой-бы объект мы ни взяли, существует тензор, компоненты 
которого в данной системе координат совпадают с компонентами
выбраного нами объекта.} 

Нужно только иметь ввиду, что как только компоненты тензора 
фиксированы в одной системе координат, то в любой другой 
системе они уже заданы <<автоматически>> в соответствии 
с формулой~(15).

\medskip

\noindent\textit{Упражнения}

\noindent\textbf{\theexer$^\ast$.} Рассмотрим тензор ранга 
(1,1), компоненты которого в некоторой системе координат 
совпадают с компонентами символа Кронекера $\delta_s^r$. 
Покажите, что в любой другой системе координат компоненты 
этого тензора также совпадают с компонентами символа Кронекера 
$\delta_s^r$. 
Это обстоятельство выражают словами <<символ Кронекера 
$\delta_s^r$ является тензором>>.
\addtocounter{exer}{1}

\noindent\textbf{\theexer$^\ast$.} Рассмотрим тензор ранга 
(2,0), компоненты которого в некоторой системе координат 
совпадают с компонентами символа Кронекера $\delta_{sr}$. 
Покажите, что в другой системе координат компоненты этого 
тензора, вообще говоря, не совпадают с компонентами символа 
Кронекера $\delta_{sr}$. 
Это обстоятельство выражают словами <<символ Кронекера 
$\delta_{sr}$ не является тензором>>.
\addtocounter{exer}{1}

\noindent\textbf{\theexer.} Выясните является ли тензором 
символа Кронекера $\delta^{sr}$.
\addtocounter{exer}{1}

\subsection*{Операции над тензорами.} 
На множестве тензоров со многими индексами имеются три 
основные операции, называемые сложением, умножением и свёрткой. 
Все они определяются путём свед\'ения к соответствующим 
операциям над объектами, определённым выше. 

Например,  пусть $x_{st}^r$ и $y_{st}^r$ --- два тензора ранга 
(2,1).
Фиксируем произвольную систему координат; 
сложим объекты, соответствующие нашим тензорам в этой системе 
координат; по получившемуся объекту построим тензор $z_{st}^r$, 
который и назовём \textit{суммой} тензоров $x_{st}^r$ и 
$y_{st}^r$.
Кратко будем писать $z_{st}^r=x_{st}^r+y_{st}^r$.

Данное выше определение суммы тензоров нуждается в проверке 
того, что результат не зависит от использованной в определении 
системы координат.
Для этого достаточно проверить, что если мы применим наше 
определение дважды (порознь в <<старой>> и <<новой>> системах 
координат), то полученные объекты $z_{st}^r$ и 
$\overline{z}_{st}^r$ будут связаны между собой тензорным 
законом преобразования (15).

Это легко проверяется непосредственным вычислением:
$$\overline{z}_{st}^r=\overline{x}_{st}^r+\overline{y}_{st}^r=
c_m^r\gamma_s^n\gamma_t^px_{np}^m+c_m^r\gamma_s^n\gamma_t^py_{np}^m=
c_m^r\gamma_s^n\gamma_t^p(x_{np}^m+y_{np}^m)=
c_m^r\gamma_s^n\gamma_t^p z_{np}^m.$$

Уравнение, связывающее компоненты тензоров, называется 
\textit{тензорным}, если из того, что оно справедливо в одной 
системе координат следует, что оно выполнено в любой системе 
координат.

Примером тензорного уравнения может служить система равенств 
$x_s^r=y_s^r$, связывающих два тензора $x_s^r$ и $y_s^r$ ранга 
(1,1).
Чтобы доказать это, достаточно заметить, что у тензора 
$x_s^r-y_s^r$ в одной из систем координат все компоненты 
равны нулю, а значит, все компоненты равны нулю и в любой 
другой системе координат.

\medskip

\noindent\textit{Упражнения}

\noindent\textbf{\theexer.}
Самостоятельно дайте определение произведения тензоров 
$x_{st}^r$ и $y^{mn}$ и свёртки тензора $x_{kst}^{rp}$ 
по индексам $k$ и $p$. 
Убедитесь, что результом будут тензоры.
\addtocounter{exer}{1}

\noindent\textbf{\theexer.}
Покажите, что $x_{st}^r y_r^p$ есть тензор ранга (2,1).
\addtocounter{exer}{1}

\noindent\textbf{\theexer$^\ast$.}
Покажите, что равенства $(a_{st}^r+b_{st}^r)x^t=d_s^r$ 
являются тензорным уравнением.
\addtocounter{exer}{1}

\subsection*{Обратный тензорный признак, или правило сокращения.} 
Мы видели, что три основные операции: сложение, умножение и 
свёртка --- являются тензорными и любая комбинация этих операций,
выполненная над тензорами, очевидно, приводит к новым тензорам.
Это позволяет распознать тензорный характер какого-нибудь 
объекта, заметив, что он образован при помощи этих операций 
над известными тензорами. 

Например, если известно, что в некоторой системе координат 
объекты $x_{st}^r$, $y^{st}$ и $z^r$ связаны соотношением
$x_{st}^ry^{st}=z^r$ и известно, что
$x_{st}^t$ и $y^{st}$ --- тензоры, то можно сразу
сказать, что $z^r$ является тензором ранга (0,1), поскольку 
он получен умножением и последующим свёртыванием двух тензоров.
Но иногда важно уметь распознавать тензоры, так сказать, 
<<обратным способом>>: а именно, если известно, что $y^{st}$ и 
$z^r$ --- тензоры, обязан ли $x_{st}^r$ быть тензором?

Сформулируем соответствующий признак, называемый <<обратным 
тензорным признаком>> или <<правилом сокращения>> для конкретного 
примера (но отметим, что он, очевидно, может быть обобщен 
на объекты любого ранга).

\textit{Пусть в каждой системе координат нам дано соотношение
$$x(r,s,t)y^{st}=z^r,\eqno(16)$$ причём известно, что 
$x(r,s,t)$ это некоторый набор чисел, занумерованных индексами 
$r$, $s$, $t$} (т.\,е. изначально нам ничего не известно о том, 
как связаны между собой наборы чисел $x(r,s,t)$ в разных системах 
координат), \textit{$y^{st}$~--- произвольный тензор ранга 
$(0,2)$ и $z^r$~--- некоторый тензор ранга} (0,1) 
(т.\,е. $z^r$ определяется однозначно как только задан тензор 
$y^{st}$.)
\textit{При сделанных предположениях $x(r,s,t)$ является 
тензором ранга} (2,1) (и, в частности, может быть записан 
в виде $x_{st}^r$).

Для доказательства запишем соотношение (16) в <<новой>> 
системе координат, а затем воспользуемся тензорным законом 
преобразования $y^{st}$ и $z^r$ и формулами преобразования (4). 
В результате получим
$$\overline{x}(r,s,t)\overline{y}^{st}=\overline{z}^r=
c_m^rz^m=c_m^rx(m,n,p) y^{np}=c_m^r x(m,n,p) \gamma_s^n\gamma_t^p \overline{y}^{st}$$
или
$$\bigl(\overline{x}(r,s,t)-c_m^r x(m,n,p) \gamma_s^n\gamma_t^p\bigr) \overline{y}^{st}=0.\eqno(17)$$

Но $y^{st}$~--- произвольный тензор, и, следовательно, 
компоненты тензора $\overline{y}^{st}$ в <<новой>> системе 
координат могут принимать любые наперёд заданные значения, 
например, такие, что именно компонета с индексами $s$ и $t$, 
равна единице, а все остальные --- нулю. 
Тогда из (17) получаем 
$$\overline{x}(r,s,t)-c_m^r x(m,n,p) \gamma_s^n\gamma_t^p=0.$$
Но это и означает, что $x(r,s,t)$ является тензором ранга (2,1)
(в частности, его следует записывать в виде~$x_{st}^r$).

\medskip

\noindent\textit{Упражнения}

\noindent\textbf{\theexer$^\ast$.}
Докажите, что если в соотношении (16) $y^{st}$ является 
симметричным, а в остаьном произвольным тензором ранга (0,2), 
то объект $x(r,s,t)+x(r,t,s)$ является тензором. 
Выведите отсюда, что если $x(r,s,t)$ есть объект, 
симметричный относительно индексов $s$ и $t$, то $x(r,s,t)$ 
есть тензор.
\addtocounter{exer}{1}

\noindent\textbf{\theexer.}
Докажите, что если в соотношении (16) $y^{st}$ является 
антисимметричным, а в остаьном произвольным тензором ранга 
(0,2), то объект $x(r,s,t)-x(r,t,s)$ является тензором. 
Выведите отсюда, что если $x(r,s,t)$ есть объект, 
антисимметричный относительно индексов $s$ и $t$, 
то $x(r,s,t)$ есть тензор.
\addtocounter{exer}{1}

\noindent\textbf{\theexer.}
Докажите, что если для любого контравариантного вектора
$x^r$ выражение $a_{rs}x^rx^s$ является скаляром
(т.\,е. тензором ранга (0,0)) и 
объект $a_{rs}$ симметричен, то $a_{rs}$ есть тензор.
\addtocounter{exer}{1}

\subsection*{Метрический тензор. Опускание и поднимание индексов.}
Из аналитической геометрии вы знаете, что скалярное 
произведение векторов $(x^1,x^2,x^3)$ и $(y^1,y^2,y^3)$ 
в прямоугольных декартовых координатах задаётся формулой
$x^1y^1+x^2y^2+x^3y^3$, 
а в косоугольной системе координат~--- формулой
$$
g_{rs}x^ry^s,\eqno(18)
$$ 
где $g_{rs}$ --- некоторые числовые коэффициенты, 
характеризующие данную косоугольную систему координат, 
причём $g_{rs}=g_{sr}$.

Поскольку скалярное произведение векторов есть инвариант 
(т.\,е. является тензором ранга (0,0)), 
то, используя правило сокращения (точнее --- упражнение 33), 
заключаем, что $g_{rs}$ есть тензор ранга (2,0).
Его называют \textit{фундаментальным} или \textit{метрическим} 
тензором (уточняя, если нужно, что речь идёт о ковариантном 
метрическом тензоре).

С помощью фундаментального тензора можно любому контравариантному 
вектору $x^r$ сопоставить ковариантный вектор $x_r$ по формуле 
$$x_r=g_{rs}x^s.\eqno(19)$$

Переход от контравариантного вектора к ковариантному, 
выполненный  помощью формулы (19), называют \textit{опусканием 
индекса.}
Очевидно, эта операция задаёт линейное отображение 
линейного пространства контравариантных векторов в
линейное пространство ковариантных векторов. 

Покажем, что \textit{оно отображает трёхмерное пространтво 
контравариантных векторов на всё трёхмерное пространство 
ковариантных векторов}, и, значит, обратимо.
Фиксировав систему координат мы можем указать три 
ковариантных вектора, координаты которых в этой системе 
координат равны $(1,0,0)$, $(0,1,0)$ и $(0,0,1)$.
Очевидно, они линейно независимы и координаты любого другого 
ковариантного вектора в этой системе координат выражаются 
через них линейно.
Следовательно, пространство ковариантных векторов трёхмерно.
Утверждение о сюръективности отображения (19) будет доказано,
если мы убедимся, что ядро отображения (19) состоит только 
из нулевого вектора.
Допустим противное, и обозначим через $x^r$ ненулевой 
контравариантный вектор из ядра, т.\,е. такой контравариантный 
вектор, что 
$$g_{rs}x^s=0.\eqno(20)$$
Умножив $r$-е равенство (20) на $r$-ю компоненту произвольного 
контравариантного вектора $y^r$ получим $g_{rs}x^sy^r=0$. 
Тем самым мы обнаружили в <<обычном>> трёхмерном пространстве 
контравариантных векторов ненулевой вектор $x^r$,
ортогональный любому вектору этого пространства. 
Пришли к противоречию, которое
доказывает обратимость отображения (19). 

Отображение, обратное к линейному, само является линейным 
и однозначно определяется  исходным отображением.
Следовательно, отображение, обратное к (19), может быть 
единственным образом записано в координатах в виде 
$$x^r=g^{rs}x_s.\eqno(21)$$

Из правила сокращения немедленно вытекает, что объект $g^{rs}$, 
определённый в координатах равенствами (21), является тензором.
Его называют \textit{фундаментальным} или \textit{метрическим} 
тензором (уточняя, если нужно, что речь идёт о контравариантном 
метрическом тензоре).
Вывод основных свойств метрического тензора $g^{rs}$ мы 
оставляем читателю в качестве упражнения (см. упражнения 35--37).

Переход от ковариантного вектора к контравариантному, 
выполненный  помощью формулы (21), называют \textit{подниманием 
индекса.}
Операции опускания и поднимания индексов естественным образом 
распространяют на тензоры произвольного ранга. 
Например, следующее равенство является просто определением:
$x_n^{mp} = g_{nr}g^{ms}g^{pt}x_{st}^r$.

\medskip

\noindent\textit{Упражнения}

\noindent\textbf{\theexer$^\ast$.}
Докажите, что равенства $g_{rs}=1$, если $r=s$ и $g_{rs}=0$, 
если $r\neq s$ являются необходимыми и достаточными для того, 
чтобы система координат была декартовой ортогональной
(т.\,е. системой, построенной по ортонормированному базису).
\addtocounter{exer}{1}

\noindent\textbf{\theexer$^\ast$.}
Докажите равенства $g_{rs}g^{st}=\delta_r^t$ и 
$g_{rs}g^{st}=\delta_r^t$, каждое из которых позволяет 
находить компоненты контравариантного метрического тензора 
$g^{st}$, если известны компоненты ковариантного 
метрического тензора $g_{rs}$.
\addtocounter{exer}{1}

\noindent\textbf{\theexer$^\ast$.}
Докажите, что если векторы $x^r$ и $x_r$ связаны между собой 
соотношениями (19) (или (21)), то имеет место равенство 
$g_{rs}x^rx^s=g^{st}x_rx_s$.
\addtocounter{exer}{1}

\noindent\textbf{\theexer$^\ast$.}
Докажите, что формула  $g^{st}x_rx_s$ задаёт скалярное 
произведение в пространстве ковариантных векторов.
Обратите внимание, что в предыдущем упражнении вы доказали, 
что отображение (19) (также как и (21)) сохраняет скалярное 
произведение (а значит, является изоморфизмом пространств 
контравариантных и ковариантных векторов).
\addtocounter{exer}{1}

\section{Вариации на тему тензоров}

Здесь мы рассмотрим понятия, родственные понятию тензора 
и находящие широкое применение в физике.
Ясно, что в определении тензора совсем 
не обязательно требовать, чтобы индексы пробегали 
значаения от 1 до 3 (как это было в предыдущем разделе). 
Разрешив им пробегать значения от 1 до некоторого 
натурального числа $n$ мы нисколько не
изменим конструкций, изложенных выше 
(даже доказательства останутся теми же).
Мы считаем такое обобщение тривиальным и не будем его 
здесь обсуждать.
Ниже мы кратко наметим подходы к более содержательным 
вариациям понятия тензора, хорошо зарекомендовавшие себя 
при решении физических задач.

\subsection*{Ортогональные тензоры.}
В физике бывают величины, которые не являются тензорами, 
но преобразуются по тензорному закону (15),
если и <<старая>> и <<новая>> системы координат являются 
ортогональными и на всех осях выбран одинаковый масштаб. 
В этом случае матрицы перехода $\gamma_r^s$ и $c_r^s$ 
от <<старого>> базиса к <<новому>> и обратно
являются ортогональными (см. упражение 21 и формулу (13)).
Объекты, преобразующиеся по тензорному закону (15), 
в котором участвуют только ортогональные матрицы
$\gamma_r^s$ и $c_r^s$, называют 
\textit{ортогональными тензорами.}

Очевидно, всякий тензор изучавшийся нами в предыдущем разделе, 
является в то же время ортогональным тензором (если запретить 
себе рассматривать его к косоугольных координатах, т.\,е. если
использовать только ортогональные матрицы $\gamma_r^s$ и $c_r^s$). 
Обратное неверно (сравните упражнения 26, 27 и 38). 

Основная особенность ортогональных тензоров состоит в том, что 
у них можно не различать верхние и нижние индексы (см. 
упражнение 39).
Именно это обстоятельство позволяло вам не различать верхние 
и нижние индексы в курсе математического анализа (там вы 
работали в арифметическом пространстве $\Bbb R^n$
с фиксрованным ортонормированным базисом).

\medskip

\noindent\textit{Упражнения}

\noindent\textbf{\theexer$^\ast$.} 
Докажите, что символы Кронекера $\delta_{rs}$ и 
$\delta^{rs}$ являются ортогональными тензорами.
\addtocounter{exer}{1}

\noindent\textbf{\theexer.} 
Пусть $x_{st}^r$ --- произвольный ортогональный тензор 
ранга (2,1) и $x_t^{rs}$ --- соответствющий ему тензор, 
полученный поднятием индекса $s$, т.\,е. 
$x_t^{rs}=g^{sm}x_{mt}^r$.
Докажите, что $x_t^{rs}=x_{st}^r$ и тем самым убедитесь, 
что для ортогонального тензора операция поднятия индекса 
не изменяет его компонет. 
Проверьте, что операция опускания индекса также не изменяет 
компонент ортогонального тензора.
\addtocounter{exer}{1}

\subsection*{Псевдотензоры.} Допустим, что 
в каждой системе координат трёхмерного евклидова пространства
задан объект $a_{i_1\dots i_m}^{j_1\dots j_n}$ ранга $(m,n)$,
т.\,е. объект, зависящий от $m$ нижних индексов и от $n$.
Предположим также, что переходе от <<старой>> сиcтемы координат 
к <<новой>>, задаваемом формулами 
$\overline{x}^r= c_s^rx^s$ (или, что то же самое, 
$x^r= \gamma_s^r\overline{x}^s$),
преобразуется по правилу
$$
\overline{a}_{i_1\dots i_m}^{j_1\dots j_n}= (\det \gamma_r^s)^M
\gamma_{i_1}^{r_1}\dots \gamma_{i_m}^{r_m}c_{s_1}^{j_1}\dots c_{s_n}^{j_n}
a_{r_1\dots r_m}^{s_1\dots s_n},$$
где $M$ --- некоторое фиксированное число.
Такой объект называют \textit{относительным тензором} 
или \textit{псевдотензором веса} $M$.

Тензоры, рассмотренные в предыдущем разделе, можно считать 
псевдотензорами веса ноль.
Если хотят отличить их от псевдотензоров, то называют 
\textit{абсолютными} или \textit{истинными} тензорами. 
Псевдотензоры первого и нулевого порядков называют 
\textit{псевдовекторами} и \textit{псевдоскалярами} 
соответственно.

На множестве псевдотензоров вводят операции сложения, 
умножения и свёртки, очень напоминающие соответствующие 
операции для истинных тензоров (см. упражнения 41--43). 
Для псевдотензоров справедлив и <<обратный тензорный признак>> 
(см. упражнение 44).

Примерами псевдотензоров могут служить символы Леви--Чивиты 
$e_{rst}$ и  $e^{rst}$, определёные выше. 
Более детально это можно описать так.
Предположим, что в каждой системе координат имеется объект, 
компоненты которого совпадают с компонентами символа 
Леви--Чивиты $e_{rst}$.
Утверждается, что при переходе от <<старой>> системы координат 
к <<новой>> компоненты этого объекта преобразуются как 
компоненты псевдотензора ранга (3,0) веса $-1$.

Чтобы доказать это, воспользуемся формулой (2). 
Мы уже доказали её для произвольного объекта $x_s^r$ ранга (1,1), 
а сейчас применим её к объекту $\gamma_s^r$, задающему
переход от <<новой>> системы координат к <<старой>>: 
$x^r=\gamma_s^r\overline{x}^s$.
В результате получим
$$
e_{rst} \gamma_m^r\gamma_n^s\gamma_p^t=e_{mnp} \det \gamma_j^k
\qquad\mbox{или}\qquad
\overline{e}_{mnp}=e_{mnp}=\bigl(\det \gamma_j^k\bigr)^{-1} \gamma_m^r\gamma_n^s\gamma_p^t e_{rst}, 
$$
что и означает, что символ Леви--Чивиты $e_{rst}$ 
преобразуется как псевдотензор веса $-1$.

\medskip

\noindent\textit{Упражнения}

\noindent\textbf{\theexer$^\ast$.} 
Докажите, что если $x$ --- псевдоскаляр веса 1 и 
$y_{st}^r$ --- псевдотензор веса~$M$,
то $x^{-M}y_{st}^r$ есть истинный тензор.
\addtocounter{exer}{1}

\noindent\textbf{\theexer.} 
Чтобы найти \textit{сумму} двух псевдотензоров одного и 
того же ранга и веса, складывают их соответствующие компоненты 
(в любой системе координат).
Докажите, что в результате получится псевдотензор того 
же ранга и веса, что и исходные слагаемые.
\addtocounter{exer}{1}

\noindent\textbf{\theexer.} 
Чтобы найти \textit{произведение} двух псевдотензоров, 
один из которых  имеет ранг $(m,n)$ и вес $M$, а другой ---  
ранг $(p,q)$ и вес $P$, умножают каждую компоненту первого 
псевдотензора на каждую компоненту второго (в любой системе 
координат).
Докажите, что в результате получится псевдотензор  
ранга $(m+p, n+q)$ и веса $M+Q$.
\addtocounter{exer}{1}

\noindent\textbf{\theexer.} 
Докажите, что в результате \textit{сворачивания} 
псевдотензора ранга $(m,n)$ и веса $M$ по верхнему и нижнему 
индексам получится псевдотензор ранга $(m-1, n-1)$ и веса $M$.
\addtocounter{exer}{1}

\noindent\textbf{\theexer$^\ast$.} 
Докажите обратный тензорный признак в следующей формулировке:
\textit{если нам дано соотношение $x(r,s,t)y^{st}=z^r$, 
причём известно, что $y^{st}$ является произвольным 
псевдотензором ранга $(0,2)$ и веса $M$, а $z^r$ --- 
определённым} (\textit{но своим для каждого} $y^{st}$) 
\textit{псевдотензором веса $N$, то объект $x(r,s,t)$ 
является псевдотензором ранга $(2,1)$ веса} $N-M$
(\textit{и тем самым его следует обозначить через} $x_{st}^r$).
\addtocounter{exer}{1}

\noindent\textbf{\theexer$^\ast$.}
Докажите, что символ Леви--Чивиты $e^{rst}$ 
является псевдотензором веса 1.\newline
{\hphantom{a}}\hfill [\texttt{Указание:} Воспользуйтесь 
упражнением 14.]
\addtocounter{exer}{1}

\noindent\textbf{\theexer.}
Покажите, что если все компоненты псевдотензора равны 
нулю в одной системе координат,
то они будут равны нулю и в любой другой системе координат.
\addtocounter{exer}{1}

\noindent\textbf{\theexer$^\ast$.}
Покажите, что если все компоненты двух псевдотензоров 
одного и того же ранга и веса
равны в одной системе координат,
то они будут равны и в любой другой системе координат.
\addtocounter{exer}{1}

\noindent\textbf{\theexer.}
Покажите, что если $x_{st}^r$ является псевдотензором 
веса $M$, то 
$x_{st}^r=\bigl(\det c_k^j\bigr)^M\gamma_m^rc_s^nc_t^p\overline{x}_{np}^m.$
\addtocounter{exer}{1}

\noindent\textbf{\theexer.} 
Докажите, что если $x_s^r$ является истинным тензором 
ранга (1,1), то $\det x_s^r$ является истинным скаляром.
\addtocounter{exer}{1}

\noindent\textbf{\theexer$^\ast$.} 
Докажите, что если $x_{rs}$ является истинным тензором 
ранга (2,0), то $\det x_{rs}$ является псевдоскаляром веса 2.
\addtocounter{exer}{1}

\noindent\textbf{\theexer.} 
Докажите, что если $x^{rs}$ является истинным тензором 
ранга (0,2), то $\det x^{rs}$ является псевдоскаляром веса $-2$.
\addtocounter{exer}{1}

\noindent\textbf{\theexer$^\ast$.} 
Докажите, что если $x_{rs}$ является истинным тензором 
ранга (2,0), причём $\det x_{kj}>0$ , то каждый из объектов 
$$
\sqrt{\det x_{kj}}e_{rst} \quad\mbox{и}\quad \frac{e^{rst}}{\sqrt{\det x_{kj}}}
$$
является истинным тензором.
\addtocounter{exer}{1}

\noindent\textbf{\theexer.} 
Предположим, что $x_{rs}$ и $y_{rs}$ являются истинными 
тензорами, $\alpha\in\Bbb R$, и
определитель $\det (x_{rs}-\alpha y_{rs})$ обращается 
в нуль при $\alpha=\alpha_0$ в некоторой системе координат. 
Докажите, что в любой другой системе координат определитель
$\det (\overline{x}_{rs}-\alpha \overline{y}_{rs})$ 
обращается в нуль при $\alpha=\alpha_0$.
Другими словами, докажите, что корни уравнения 
$$\det (x_{rs}-\alpha y_{rs})=0$$ являются инвариантами.
\addtocounter{exer}{1}

\subsection*{Тензор Леви--Чивиты и векторное произведение.}
Пусть $g_{rs}$ является метрическим тензором.
Положим $g=\det g_{rs}$. 
Согласно упражнению 50, $g$ является псевдоскаляром веса 2.
Значит, выбрав в качестве <<новой>> системы координат 
декартову ортогональную (в которой, 
$\overline{g}_{rs}=\delta_{rs}$, а значит, $\overline{g}=1$),
можем написать $1=\overline{g}=\bigl(\det \gamma_s^r\bigr)^2 g$.
Отсюда следует, что $g>0$, ведь матрица $\gamma_r^s$, задающая, 
в соответствии с (5), переход от <<новой>> системы  координат 
к <<старой>>, невырождена, а значит $\det \gamma_s^r\neq 0$.

В соответствии с упражнением 52, объекты
$$
\varepsilon_{rst}=\sqrt{g}e_{rst} \quad\mbox{и}\quad 
\varepsilon^{rst}=\frac{e^{rst}}{\sqrt{g}}
$$
являются истиными тензорами. 
В литературе каждый из этих тензоров называют или
\textit{тензором Леви--Чивиты}, или
\textit{абсолютно антисимметричным единичным тензором}, или
\textit{полностью антисимметричным единичным тензором}.
Мы будет придерживаться первого названия.

Имея два контравариантных вектора $x^r$ и $y^r$ построим 
новый контравариантный вектор $z^r$ по формуле 
$$
z^r=\varepsilon^{rmn}g_{ms}g_{nt}x^sy^t.\eqno(22)
$$
Убедимся, что \textit{вектор $z^r$, построенный по формуле} (22), 
\textit{является векторным произведением векторов $x^r$ и $y^r$.}
Поскольку (22) является тензорным равенством, его достаточно 
проверить в какой-нибудь одной системе координат. 
Сделаем это в декартовой ортогональной системе с равными 
масштабами по осям. В такой системе координат 
$g_{sm}=\delta_{sm}$ и можно не различать верхние и нижние 
индексы. 
Поэтому равенство (22) примет вид $z^r=e^{rmn}x^my^n$ или, 
в развёрнутом виде,
\begin{align*}
z^1&=x^2y^3-x^3y^2,\\
z^2&=x^3y^1-x^1y^3,\tag{23}\\
z^3&=x^1y^2-x^2y^1.
\end{align*}
Теперь сразу видно, что формулы (23) являются покоординатной 
записью хорошо известной формулы 
$$
\boldsymbol{z}=
\det
\begin{pmatrix}
\boldsymbol{e}_1 &\boldsymbol{e}_2 &\boldsymbol{e}_3\\
x^1 & x^2 & x^3\\
y^1 & y^2 & y^3
\end{pmatrix}\eqno(24)
$$
для нахождения векторного произведения 
$\boldsymbol{z}=\boldsymbol{x}\times \boldsymbol{y}$ 
в системе координат, соответствующей единичному 
правому ортонормированному координатному базису 
$\boldsymbol{e}_1$, $\boldsymbol{e}_2$, $\boldsymbol{e}_3$. 
Тем самым формула (22) доказана. 

Преимущество формулы (22) по сравнению с более привычной 
формулой (24) состоит в том, что формула (22) годится 
для любой системы координат (в том числе косоугольной и левой).

Используя формулы (18) и (22) для скалярного и векторного 
произведения векторов, немендленно получаем, что формула
$$
(\boldsymbol{x},\boldsymbol{y},\boldsymbol{z})=
\boldsymbol{x}\cdot(\boldsymbol{y}\times \boldsymbol{z})=
\varepsilon^{mnp}g_{mr}g_{ns}g_{pt}x^ry^sz^t
$$
задаёт \textit{смешанное} произведение векторов 
$x^r$, $y^r$ и $z^r$ в произвольной системе координат.

\medskip

\noindent\textit{Упражнения}

\noindent\textbf{\theexer$^\ast$.} 
Докажите, что тензор $\varepsilon^{rst}$ получается из 
тензора $\varepsilon_{rst}$ поднятием индексов.

[\texttt{Указание:} заметьете, 
что вас просят доказать равенство 
$\varepsilon^{rst}=\varepsilon_{mnp}g^{rm}g^{sn}g^{tp}$, 
причём сделать это достаточно хотя бы в одной системе координат; 
проверьте его в декартовой ортогональной системе координат
с одинаковыми масштабами по осям (т.\,е. такой, что 
$g_{rm}=\delta_{rm}$).]
\addtocounter{exer}{1}

\noindent\textbf{\theexer.} 
Докажите, что компоненты тензора Леви--Чивиты 
$\varepsilon_{rst}$ в любой системе координат
могут быть получены как смешанное произведение 
$(\boldsymbol{e}_r, \boldsymbol{e}_s, \boldsymbol{e}_t)$ 
векторов базиса $\boldsymbol{e}_1$, $\boldsymbol{e}_2$, 
$\boldsymbol{e}_3$, задающего эту систему координат, 
т.\,е. докажите равенство 
$\varepsilon_{rst}=(\boldsymbol{e}_r, \boldsymbol{e}_s, \boldsymbol{e}_t)$.
\addtocounter{exer}{1}

\noindent\textbf{\theexer.}
Используя формулу (22) и упражнение  16
докажите известное правило раскрытия двойного векторного 
произведения:
$\boldsymbol{x}\times(\boldsymbol{y}\times\boldsymbol{z})=
\boldsymbol{y}(\boldsymbol{x},\boldsymbol{z})-\boldsymbol{z}(\boldsymbol{x}, \boldsymbol{y})$.
\addtocounter{exer}{1}

\subsection*{4-векторы.}
В специальной теории относительности событие задают 
четырьмя числами: тремя координатами $x$, $y$, $z$, 
описывающими положение той точки трёхмерного
евклидова пространства, где это событие произошло, 
и моментом времени $t$, соответствующим этому событию.
Эти четыре числа записывают в виде вектора 
$$(x^0,x^1,x^2,x^3)=(ct,x,y,z),$$
где множитель $c$ обозначает скорость света и введён 
для того, чтобы все координаты имели одну и ту же размерность 
(а именно, размерность длины).
Сумму векторов и произведение вектора на вещественное 
число определяют покомпонентно, а именно вектор 
$$(x^0+y^0,x^1+y^1,x^2+y^2,x^3+y^3)$$ называют 
\textit{суммой} векторов 
$$\boldsymbol{x}=(x^0,x^1,x^2,x^3)
\qquad\mbox{и}\qquad 
\boldsymbol{y}=(y^0,y^1,y^2,y^3),$$ 
а вектор 
$(\alpha x^0,\alpha x^1,\alpha x^2,\alpha x^3)$ 
называют результатом умножения вектора $\boldsymbol{x}$ 
на вещественное число $\alpha$.
Легко видеть, что совокупность векторов $\boldsymbol{x}$ 
с так определёнными операциями сложения и умножения на числа 
является линейным (или, что то же самое, векторным) пространством. 

Каждым двум векторами $\boldsymbol{x}$ и $\boldsymbol{y}$
этого пространства сопоставляют число 
$\boldsymbol{x}\cdot\boldsymbol{y}$, заданное формулой
$$
\boldsymbol{x}\cdot\boldsymbol{y}=x^0y^0-x^1y^1-x^2y^2-x^3y^3.\eqno(25)
$$
 
Отображение, заданное формулой (25), является линейным по
первому аргументу и симметричным, но не является положительно 
определённым (зато оно обладает следующим свойством 
невырожденности:
\textit{если $\boldsymbol{x}\cdot\boldsymbol{y}=\boldsymbol{0}$ для всех $\boldsymbol{y}$,
то} $\boldsymbol{x}=\boldsymbol{0}$).
Тем самым для $\boldsymbol{x}\cdot\boldsymbol{y}$ выполнены две
аксиомы скалярного произведения, а вместо третьей выполняется 
аксиома невырожденности.
Такую операцию, сопоставляющую паре векторов $\boldsymbol{x}$ 
и $\boldsymbol{y}$ число $\boldsymbol{x}\cdot\boldsymbol{y}$ 
по формуле (25), называют \textit{псевдоскалярным произведением} 
векторов.

Описанное выше пространство векторов 
$\boldsymbol{x}=(x^0,x^1,x^2,x^3)$
с псевдоскалярным произведением (25) называют 
\textit{пространством Минковского} или 
\textit{псевдоевклидовым пространством}. 
Мы будем придерживаться первого названия, отдавая честь 
Г\'ерману Минк\'овскому (Hermann Minkowski; 1864--1909) 
--- немецкому математику, с успехом применявшему 
геометрические методы для решения трудных задач в самых разных 
разделах математики и физики, включая теорию чисел, 
математическую физику и специальную теорию относительности.

До сих пор наши построения происходили в какой-то 
фиксированной системе координат, порождённой некоторым базисом 
$\boldsymbol{e}_0$, $\boldsymbol{e}_1$, $\boldsymbol{e}_2$, 
$\boldsymbol{e}_3$, который называют \textit{стандартным}. 
В частности, составляющие его векторы попарно ортогональны и
$$
-\boldsymbol{e}_0\cdot\boldsymbol{e}_0=
\boldsymbol{e}_1\cdot\boldsymbol{e}_1=
\boldsymbol{e}_2\cdot\boldsymbol{e}_2=
\boldsymbol{e}_3\cdot\boldsymbol{e}_3=-1.
$$

Любая система координат, порождённая базисом 
$\overline{\boldsymbol{e}}_0$, $\overline{\boldsymbol{e}}_1$, 
$\overline{\boldsymbol{e}}_2$, $\overline{\boldsymbol{e}}_3$, 
таким, что составляющие его векторы попарно ортогональны и
$-\overline{\boldsymbol{e}}_0\cdot\overline{\boldsymbol{e}}_0=
\overline{\boldsymbol{e}}_1\cdot\overline{\boldsymbol{e}}_1=
\overline{\boldsymbol{e}}_2\cdot\overline{\boldsymbol{e}}_2=
\overline{\boldsymbol{e}}_3\cdot\overline{\boldsymbol{e}}_3=-1$,
называется \textit{галилеевой}\footnote{Галил{\'е}о Галил{\'е}й 
(Galileo Galilei; 1564 -- 1642) --- итальянский физик, 
механик, астроном, философ и математик, оказавший 
значительное влияние на науку своего времени.}
 системой координат.
Роль галилеевых систем координат в псевдоевклидовом 
пространстве аналогична роли декартовых ортогональных 
систем координат в евклидовом пространстве.
В специальной теории относительности инерциальная 
система отсчёта обычно задаётся галилеевой системой координат.

При переходе от одной галилеевой системы координат 
$(x^0,x^1,x^2,x^3)$ 
к другой галилеевой системе координат 
$(\overline{x}^0,\overline{x}^1,\overline{x}^2,\overline{x}^3)$ 
координаты пересчитываются посредством некоторой матрицы перехода
$$
\overline{x}^r= c_s^rx^s.\eqno(26)
$$
Эти матрицы перехода $c_s^r$ во многом аналогичны 
ортогональным матрицам.
Главное в этой аналогии в том, что матрицы перехода сохраняют 
псевдоскалярное произведение (25) и образуют группу 
(т.\,е. произведение двух  матриц перехода снова является
матрицей перехода). Эту группу называют \textit{группой Лоренца}
в честь выдающегося голландского физика Х\'ендрика Л\'оренца 
(Hendrik Lorentz; 1853--1928).

Теперь мы можем сформулировать определение, являющееся 
главной целью этого раздела: 
говорят, что в пространстве Минковского задан 4-\textit{вектор}, 
если в каждой галилеевой системе координат задан объект $x^r$, 
компоненты которого при переходе к любой другой
галилеевой системе координат преобразуются по формулам~(26).

Другими словами можно сказать, что 4-векторы аналогичны 
ортогональным контравариантным
тензорам ранга (0,1) в евклидовом пространстве. 
Безусловно, можно рассматривать подобные объекты с любым 
числом верхних и нижних индексов. Но при изучении 
электродинамики и специальной теории относительности
вы познакомитесь с понятиями 4-скорости, 4-ускорения, 4-момента, 
4-силы, 4-тока, и 4-потенциала и узнаете, что все они 
являются 4-векторами (см.,  например, \cite{GP10}).

\medskip

\noindent\textit{Упражнения}

\noindent\textbf{\theexer.} 
Пусть $c_s^r$ есть матрица перехода от одной галилеевой 
системы координат 
к другой галилеевой системе координат. Докажите, что 

$$
c_s^0c_r^0-c_s^1c_r^1-c_s^2c_r^2-c_s^3c_r^3=
\begin{cases}
0,&\mbox{если \ } s\neq r;\\
1,&\mbox{если \ } s=r=0;\\
-1,&\mbox{если \ } s=r\neq 0.
\end{cases}
\eqno(27)
$$
\addtocounter{exer}{1}

\noindent\textbf{\theexer.} 
Пусть даны две системы координат, одна из которых является 
галилеевой, а вторая --- произвольной. И пусть известно, 
что матрица перехода $c_s^r$  от первой системы
координат ко второй удовлетворяет соотношениям (27). 
Докажите, что тогда вторая система координат с необходимостью 
является галилеевой.
\addtocounter{exer}{1}

\noindent\textbf{\theexer.}
В специальной теории относительности показывается, что 
при переходе от одной галилеевой системы координат 
$(x^0,x^1,x^2,x^3)$ 
к другой галилеевой системе координат 
$(\overline{x}^0,\overline{x}^1,\overline{x}^2,\overline{x}^3)$, 
движущейся относительно первой со скоростью $v$ параллельно 
оси $x^1$,  координаты пересчитываются по формулам
$$
\overline{x}^0=\frac{x^0-(v/c)x^1}{\sqrt{1-(v/c)^2}},\quad
\overline{x}^1=\frac{x^1-(v/c)x^0}{\sqrt{1-(v/c)^2}},\quad
\overline{x}^2=x^2,\quad
\overline{x}^3=x^3.
$$
Убедитесь, что матричная запись (26) этого преобразования координат
осуществляется с помощью следующей матрицы $c_s^r$:
$$
\begin{pmatrix}
\frac{1}{\sqrt{1-(v/c)^2}}    & -\frac{v/c}{\sqrt{1-(v/c)^2}} & 0 & 0\\
-\frac{v/c}{\sqrt{1-(v/c)^2}} & \frac{1}{\sqrt{1-(v/c)^2}}    & 0 & 0\\
0                             & 0                             & 1 & 0\\
0                             & 0                             & 0 & 1
\end{pmatrix}.\eqno(28)
$$
Проверьте также, что для матрицы (28) выполняются условия (27).
\addtocounter{exer}{1}

\noindent\textbf{\theexer.}
В условиях предыдущей задачи докажите, что существует 
единственное вещественное число
$\psi$ такое, что выполнены сразу два равенства
$$
\mbox{sh\,}\psi =\frac{v/c}{\sqrt{1-(v/c)^2}} \quad\mbox{и}\quad
\mbox{ch\,}\psi =\frac{1}{\sqrt{1-(v/c)^2}}.
$$
Величину $\psi$ называют \textit{углом гиперболического 
поворота}.
Проверьте, что с её помощью матрица (28) записывается 
более компактно:
$$
\begin{pmatrix}
\mbox{ch\,}\psi   & \mbox{sh\,}\psi  & 0 & 0\\
\mbox{sh\,}\psi   & \mbox{ch\,}\psi  & 0 & 0\\
0                 & 0                & 1 & 0\\
0                 & 0                & 0 & 1
\end{pmatrix}.
$$
\addtocounter{exer}{1}

\section{Заключение}

Читателю, желающему более глубоко освоить тензорное исчисление, мы рекомендуем 
обратиться к чтению <<толстых учебников>>. 
Электронные варианты многих (может быть, даже всех) стандартных учебников по теории тензоров 
есть, например,  на сайте \texttt{http://eqworld.ipmnet.ru/ru/library/mathematics/difgeometry.htm}.

Наше изложение наиболее близко к принятому в классической книге А.\,Дж.~Мак-Коннела \cite{MC63}.
В его основе лежит так называемый координатный подход к тензорной алгебре и тензорному анализу.
Он компактен, гибок и широко используется в физической литературе.
Однако существует и другой, бескоординатный, подход к понятию тензора, с которым 
можно познакомиться, например, по книге А.\,И.~Кострикина и Ю.\,И.~Манина \cite{KM80}.
Бескоординатному подходу отдаётся предпочтение в некоторых продвинутых областях современной
теоретической физики.

Тем, кто не готов сразу читать <<толстые учебники>>, 
мы рекомендуем пособие В.\,А.~Топоногова
\cite{To95}\footnote{Это пособие опубликовано также в Интернете:
\texttt{http://www.phys.nsu.ru/ok03/manuals\_1995-1996.html}}, 
написанное 15 лет назад специально для студентов физического факультета 
Новосибирского государственного университета.

Пользуясь случаем автор выражает благодарность Д.\,Ю.~Иванову, 
стимулировавшему его к написанию этой статьи и оказавшему 
существенное влияние на отбор материала.

\bigskip

{\small
\hfill\textit{Материал поступил в редколлегию 04.12.2011}

\vskip2cm

\textbf{Victor Alexandrov}

\medskip

\textbf{THE FIRST ACQUAINTANCE WITH THE TENSOR}

\bigskip

The main concepts of the theory of tensors are presented.
The emphasis is on the basic notions of tensor algebra and
practical skills in culculations involving the Kronecker 
delta and Levi-Civita symbol. 
Sixty routine exercises are included.
The article is intended for junior students of physical,
mathematical, and geophysical departments of universities.

\textit{Keywords:} tensor, Kronecker delta, Levi-Civita symbol.
}%endsmall

\end{document}